\documentclass[12pt]{amsart}
\usepackage{amssymb}
\newtheorem{theorem}{Theorem}

\newtheorem{corollary}[theorem]{Corollary}
\newtheorem{proposition}[theorem]{Proposition}
\newtheorem{lemma}[theorem]{Lemma}
\theoremstyle{definition}
\newtheorem{conjecture}[theorem]{Conjecture}
 \theoremstyle{remark}
\newtheorem{example}[theorem]{\rm{EXAMPLE}}
\newtheorem{definition}[theorem]{\rm{DEFINITION}}
\newtheorem{remark}[theorem]{\rm{REMARK}}

\def\Irr{{\rm Irr}}
\def\Bl{{\mbox{\rm Bl}}}

\def\k{{\rm k}}
\def\Cl{{\rm Cl}}
\def\k{{\rm k}}
\def\IBr{{\mbox{\rm IBr}}}

\def\l{{\rm l}}
\begin{document}

\title{Real characters in blocks}

\author{L.~H\'ethelyi}
\address{Department of Algebra, Budapest University of Technology and Economics, H-1521 Budapest, M\H uegyetem rkp.~3--9.}
\email{hethelyi@math.bme.hu}

\author{E.~Horv\'ath}
\address{Department of Algebra, Budapest University of Technology and Economics, H-1521 Budapest, M\H uegyetem rkp.~3--9.}
\email{he@math.bme.hu}

\author{E.~Szab\'o}
\address{Alfr\'ed R\'enyi Mathematical Institute of the Hungarian Academy of Sciences,
Budapest, Re\'altanoda u. 13-15.}
\email{endre@renyi.hu} 
\abstract{%
We consider  real versions of  Brauer's $\k(B)$ conjecture, Olsson's conjecture and
Eaton's conjecture. We prove the real version of Eaton's  conjecture
for $2$-blocks of groups with cyclic defect group 
and for  
the principal $2$-blocks of groups with trivial  real core. 
We also characterize $G$-classes,
 real and rational  $G$-classes  
 of the defect group of $B$.
}}

\maketitle
\footnotetext[1]{2010 Mathematics Subject Classification, 20C20, 20C15}
\section{Introduction}
Several authors have been investigating real classes, characters
and blocks of finite groups, see e.g. 
\cite{Brau,Bra,DolNavTie,Gow,G-M,HH,IMN,IN,Iwasaki,K-N,MN,M,Mur,NavSanTie, NavSanTie2}. 
The aim of this note is to formulate  real versions of Brauer's
$\k(B)$ conjecture, see \cite{Bra56}, Olsson's conjecture, see \cite{Ol84},
 and Eaton's conjecture, see \cite{Eat03}, for
$2$-blocks.
   We
 give special cases when we can prove the real versions of  them. 
The last part of the paper deals with fusion of elements
 of  defect groups.

\section{Notations and terminology}

In this note $G$ will always denote a finite group, $p$ a prime integer, which is $2$
except for the last section of  the paper.
Let $(R,k,F)$ be a $p$-modular system, where $R$ is a complete discrete valuation
ring with quotient field $k$ of characteristic zero and residue class field $F$
 of characteristic $p$. We assume that $k$ and $F$ are splitting fields of all the
subgroups of $G$. We may also assume that $k$ is a subfield of the
complex numbers.
Complex conjugation  acts on $\Irr(G)$. A character is {\bf real} if it is conjugate
to itself, in other words if it is real valued.
An element of $G$ is  {\bf real} if it is conjugate to its inverse.
An element $x$ of a subgroup $H$ of $G$ we call {\bf $H$-real}, if it can be conjugated to its inverse
by an element from the subgroup $H$.

We say that the {\bf  conjugacy class $C$ is  real } if it is equal to the class of
the inverse elements of the class. We use the notation $\Cl_{r}(G)$ for  the set of these
classes.
A {\bf $p$-block  $B$ }is  called {\bf real } if it contains the complex conjugate of
an irreducible ordinary character (and hence of all irreducible characters) 
in the block. It is known,  see e.g. \cite[Thm.~3.33]{Nav}, that a real
$2$-block always contains real valued irreducible ordinary and Brauer characters, as well.
We use the notation $\Irr_{rv}(G)$ and $\Irr_{rv}(B)$ for the set of real valued irreducible ordinary
characters in $G$ and in $B$, respectively. Let $\k_{rv}(G)$ and $\k_{rv}(B)$
stand for the sizes of these sets.
We use the notation ${\k}_{{i},{rv}}(B)$  for the number of
real valued irreducible characters of height $i$ in the $p$-block $B$. 
By Brauer's permutation lemma the number of real conjugacy classes of the group
 $G$ is equal to $k_{rv}(G)$. We use  the notation $\Bl(G|D)$ for the set of $p$-blocks
of
$G$ with defect group $D$, $D^{(n)}$ stands for the $n$-th  derived subgroup
of $D$. For constructing examples we used  the GAP system, see \cite{GAP}, and
we also  describe these groups with their GAP notation.

\section {The real  conjectures}

Unless otherwise stated, let $p=2$. Let $G$ be a finite group,  $B$  a real
$2$-block of $G$ with defect group $D$.

\begin{conjecture}[Weak real version of Brauer's conjecture]\label{conj0}
We conjecture that 
 $\k_{rv}(B)$ 
is bounded from above by the
number of $G$-real elements of  $D$.
\end{conjecture}

\begin{conjecture}[Strong real version of Brauer's conjecture]\label{conj1}
 We conjecture that  $\k_{rv}(B)$ 
 is bounded from above by the 
number of $N_G(D)$-real elements of  $D$.
\end{conjecture}

\begin{conjecture}[Real version of Olsson's conjecture]\label{conj2}
 
We conjecture that ${\k_{{0},{rv}}}(B)$ 
 is bounded from above by
the number of $N_G(D)/D'$-real elements of   $D/D'$.
\end{conjecture}

\begin{conjecture} [Real version of Eaton's conjecture]\label{conj3}

We conjecture that   
$\sum_{i=0}^n {\k}_{{i},{rv}}(B)$ 
is bounded from above by the number of

 $N_G(D)/D^{(n+1)}$-real elements 
 of  $D/D^{(n+1)}$.
\end{conjecture}

\begin{remark} One could not replace in  Conjecture \ref{conj1} the $N_G(D)$ by $D$.
The smallest example is a group of order $24$ which is the pullback of maps
$S_3\rightarrow C_2$ and $Q_8\rightarrow C_2$, (with GAP notations it is 
 $\texttt{SmallGroup}(24,4)$). In this group  there are two $2$-blocks. The nonprincipal
block $B$ has a normal defect group $D\simeq C_4$, where there are just two 
$D$-real elements, however $\k_{rv}(B)=4$. (In fact in this group all characters in $\Irr(G)$ are real). 
However, we do not know any such example for the principal block,
 or for blocks of maximal defect.
\end{remark}

\begin{remark}
If every irreducible character is real in the group $G$ then we get   
 stronger versions of the  non-real conjectures, see Remark \ref{stronger}, namely
$\k(B)$ (${\k_0}(B)$, $\sum_{i=0}^{n}{\k_i}(B))$  are bounded from above  by the number of
 elements of the defect group $D$ of $B$,($D/D'$ or $D/D^{(n+1)}$) that are
real inside $N_G(D)$ ($N_G(D)/D'$ or $N_G(D)/D^{(n+1)}$) respectively.
Of course Conjecture~\ref{conj3} implies Conjectures  \ref{conj0}, \ref{conj1}
and \ref{conj2}. 
\end{remark} 

\begin{remark}\label{stronger}
If every irreducible character of a group $G$ is real, it does not follow that
 the normalizer of its Sylow $2$-subgroup also has this property.
Let $G=\texttt{SmallGroup}(96,185)$. This group has selfnormalizing Sylow $2$-subgroups.  In the principal block of $G$ all  the $14$
 irreducible characters are real valued, its  Sylow $2$-subgroup has also $14$ irreducible characters, but only $12$ of them are real. 
This example also shows that an element can
be real in one of the Sylow $2$-subgroups,
but not real in an other Sylow $2$-subgroup,
since one can find order $4$ elements in this  with that property.
In this group all $32$ elements of the Sylow $2$-subgroup are real in $G$, but
only $28$ of them are real in $N_G(S)=S$.
\end{remark}

In the next remark we show that the $p$-analogue of Conjecture \ref{conj0}
is not true for $p>2$, and since  the defect group is abelian the other Conjectures \ref{conj1}, \ref{conj2} and \ref{conj3} also cannot hold:
 
\begin{remark}\label{forp}
Let $G=M_{11}$, $p=11$ and let $B$ be the principal block.
 Then $|D|=11$, $\k_{rv}(B)=3$,  but in $D$ there is  only one 
$G$-real element. This group is also an example for the fact that the number of real
valued irreducible characters can be different in  the Brauer correspondent blocks
with cyclic defect group if $p>2$.
 Let $b\in Bl(N_G(D)|D)$ be the Brauer correspondent of $B$. Then $\k_{rv}(b)=1$. 
If $p=2$ and the defect group is noncyclic then Brauer correspondent blocks might have
differrent number of real valued irreducible characters: let us take the same
 group $G$, then the principal $2$-block has $6$, however
 its Brauer correspondent block has $5$ real valued 
irreducible characters.
\end{remark}

\section{Nilpotent groups, symmetric groups and blocks with central defect groups}

\begin{proposition}[The nilpotent groups]\label{nilp} A stronger form of Conjecture~\ref{conj3} (hence \ref{conj2}, \ref{conj1} and \ref{conj0})
holds for nilpotent groups. If $G$ is either a $2$-group or abelian, then in Conjecture
\ref{conj2} there is equality.
\end{proposition}

\noindent
Proof.
If $G$ is nilpotent then every $2$-block  is of maximal defect, and by \cite{Gow}
 the only real $2$-block of maximal defect is the principal block $B_0$. Then $\Irr(B_0)=\Irr(G_2)$,
where $G_2\in Syl_2(G)$.
Characters of height $n$ of $G_2$ are those of degree $2^n$. This is at most the $n$th
character degree of $G_2$. By \cite[Lemma~5.12]{Is}, all irreducible characters
of height at most $n$ contain 
${G_2}^{(n+1)}$ in their kernels, hence
 $\sum_{i=0}^n {\k}_{{i},{rv}}(B_0)\leq
|\Irr_{rv}(G_2/{{G_2}^{(n+1)}})|$, which is at most the number of $G_2/G_2^{(n+1)}$-real elements
in $G_2/{{G_2}^{(n+1)}}$.

\begin{proposition}[The symmetric groups]\label{sym}
Conjectures \ref{conj1} and \ref{conj2} hold for the symmetric groups.
\end{proposition}

\noindent
Proof.
\begin{itemize}
\item[(a)]Since every irreducible character of the symmetric group is real valued and
since its Sylow $2$-subgroup  also  has this property by \cite[Thm.~4.4.8]{J-K},
Conjecture~\ref{conj1} for the principal $2$-blocks 
 reduces in this case to the non-real $\k(B)$ conjecture, which  
holds by \cite{Ol84}.
In \cite{Ol76} it is proved that the defect group $D$ of each block $B$ of weight
$w$ of $S_n$ is isomorphic to the Sylow $p$-subgroup of $S_{pw}$ and there is a canonical
height preserving bijection between the irreducible characters of $B$ and that
of the principal block of $S_{pw}$. Thus  if $p=2$ then in $D$   each element is also real, 
and again by \cite{Ol84}, Conjecture~\ref{conj1} holds also for nonprincipal
$2$-blocks of $S_n$.
\item[(b)] 
Olsson's conjecture also holds  for $S_n$ by \cite{Ol84}. Thus by similar arguments as above,    Conjecture~\ref{conj2}, also holds.
\end{itemize}

\begin{remark}
A positive answer to Eaton's conjecture for $S_n$, would imply a positive answer to
Conjecture \ref{conj3}.
\end{remark}

\begin{proposition}[Blocks with central defect groups]\label{central}
Conjecture \ref{conj3} (and hence, Conjectures \ref{conj0},\ref{conj1} and \ref{conj2}) holds for
central defect groups.
In fact we prove a slightly stronger statement:
  the strong forms of the conjectures holds for $2$-blocks with defect group $D$, where $G=DC_G(D)$:
\end{proposition}

\noindent
Proof.
Let $B\in \Bl(G|D)$ be a $2$-block of $G$, where $G=DC_G(D)$.
By \cite[Thm.~9.12]{Nav} $|\IBr(B)|=1$ and there is a bijection between
$\Irr(D)$ and $\Irr(B)$ mapping $\zeta $ to $\theta_{\zeta}$, where
$\theta_{\zeta}(g)=\zeta(g_2)\theta(g_{2'})$, if $g_2\in D$, otherwise it is zero.
 Here $\theta $ is the unique
 character in $\Irr(B)$  containing $D$ in its kernel, and $\IBr(B)=\{\theta^{0}\}$.
Moreover $ht(\theta_{\zeta})=n$ iff $\zeta(1)=2^n$.
If $B$ is a real $2$-block then $\theta $ is a real valued character and $\theta_{\zeta}$
is real valued if and only if $\zeta $ is real valued.
Thus $\k_{rv}(B)=\k_{rv}(D)$ and $\sum_{i=0}^n k_{i,rv}(B)=\sum_{i=0}^n
k_{i,rv}(D)\leq |\Irr_{rv}(D/D^{(n+1)}|$ by Proposition \ref{nilp}, this is at most
the number of $D/D^{(n+1)}$-real elements of $D/D^{(n+1)}$.

\begin{remark}
It is easy to see that if  the above conjectures  are true for the direct factors of a group  then  they are  also true for the 
  direct product:  a tensor product of characters is real iff each component is
real, if we have  defect classes $C_1\in \Cl(B_1,D_1)$ and $C_2\in \Cl(B_2,D_2)$,
 then the pair $(c_1,c_2)\in C_1\times C_2$ is a defect class of $B_1\otimes B_2$.
The defect of the character $\chi_1\otimes \chi_2$ is the sum of the defects of $\chi_1$ and $\chi_2$. The height of the product of characters is the sum of the
heights. The number of real elements in $D_1\times D_2$ is just the product
of the numbers of real elements in the direct components. 
\end{remark}

\section{Blocks with cyclic defect groups}

For a block $B\in\Bl(G)$ we consider the pairs $(x,\theta)$ with $x\in G$
a $p$-element $\theta\in\IBr(b)$, where $b\in \Bl(C_G(x))$ such that
 $b^G=B$.
As in \cite{Mur}, we call the $G$-conjugacy classes of these pairs,
denoted by $(x,\theta)^G$, the {\bf columns of $B$}. 
A  {\bf column}  $(x,\theta)^G$ is called ${\bf real}$ if $(x,\theta)^G=(x^{-1},\overline{\theta})^G$.
In \cite[Lemma 1.1]{Mur} it is proved that $\k_{rv}(B)$ is equal to the number of real
columns of $B$.

 We will  use  Dade's description \cite[Thm 68.1]{Dor} of $p$-blocks with 
cyclic defect
 groups only  for the special case $p=2$: 

Let $B$ be a $2$-block with cyclic defect group $D=\langle x \rangle$ of order
 $2^a$, $D_i=\langle x^{2^i}\rangle$, $C_i=C_G(D_i),N_i=N_G(D_i)$, for 
$i=0,\cdots,a-1$.  Let $B_0\in \Bl(N_G(D)|D)$ be the Brauer correspondent block of $B$. Let $b_0\in \Bl(C_G(D)|D)$ with $b_0^{N_0}=B_0$.
Such blocks are conjugate  in $N_0$.
 Similarly let $b_i={b_0}^{C_i}$, then every block of $C_i$ that induces $B$
 is conjugate to $b_i$ in $N_i$.
 Let $\theta^i$ be the unique irreducible Brauer character of $b_i$ for $i=0,\cdots,a-1$. 
The inertia subgroup of $b_i$ in $N_i$ is  $C_i$
for $i=0,\cdots,a-1$,  and also $|\IBr(B)|=1$. Let $\IBr(B)=\{\phi\}$.
\bigskip
 First we prove: 

\begin{lemma}\label{cyclicdefect}
 With the notation above, we have that  
 $(x,\theta^0)$, $(x^k,\theta^0)$ for $k$ odd,
$(x^2,\theta ^1)$, $(x^{2k}, \theta^1)$ for $k$ odd,
...
$(x^{2^{a-1}},\theta^{a-1})$,
$(x^{2^{(a-1)}k},\theta^{a-1})$ for $k$ odd,
  and $(1, \phi)$ are representatives of the columns of $B$.
If $x^{2^i}$ is the smallest power of $x$ which is $G$-real then  representatives of the real columns of $B$ are among
those columns whose first component is a power of $x^{2^i}$. 
\end{lemma}

\noindent
Proof.
If the first components of two pairs  generate  different subgroups, then they cannot be conjugate.
Let us take the pair $(y,\psi)$, where $y$ generates $D_j$ and  $j< a$. Then the block of $\psi$
is conjugate to $b_j$ in $N_j$, so $\psi$ is conjugate to $\theta^j$ in $N_j$.
The conjugation takes $y$ to another generator of $D_j$, i.e. to $x^{{2^j}k}$, where
$k$ is odd.
If the first component is $1$, then the second component must be $\phi $.

\begin{corollary}
Let $G$ be a finite group, let $B$ be a real $2$-block with cyclic defect group $D$.
Then Conjecture \ref{conj0} holds for $G$.
\end{corollary}

\noindent
Proof.
 We use \cite[Lemma 1.1]{Mur}, Lemma \ref{cyclicdefect} and the notations above.
Then the number of $G$-real elements in $D$ is exactly $2^{a-i}$.

We have
 that the representatives of  real columns of $B$
are  $(1,\phi)$ and  some of those columns whose first component is an element of
 $D_i$ and if it generates   $D_j$  then 
 the second component is  $\theta^j$.  Their number is at most the number of elements of $D_i$, which is
  $2^{a-i}$.

\begin{corollary} \label{normal}
Let $D$ be a cyclic normal $2$-subgroup of $G$.
Then Conjecture \ref{conj1} and hence Conjecture \ref{conj3} also
holds for blocks $B\in Bl(G|D)$.
\end{corollary}

\begin{remark}
Using similar arguments for the $p>2$ case, one gets for block with cyclic defect groups that $\k_{rv}(B)\leq \l(B)\cdot
|\{$ $G$-real elements in $D$ $\}|$.
This could be considered as some kind
 of real analogue of the so called ``Trace inequality'' in \cite[Prop. 2, p. 272]{Ol80}.
\end{remark}
\bigskip

To prove  Conjecture \ref{conj1} for $2$-blocks with cyclic defect group 
 we will need the following lemma (the $p$-analogue of it 
for $p>2$ is not true, and if $p=2$, but the defect group is noncyclic 
then the analogous result is not true either, see Example \ref{forp}):

\begin{lemma}\label{samenumber}
Let $G$ be a finite group, let  $B\in \Bl(G|D)$ be a  real $2$-block with cyclic
defect group $D$ and let $B_0\in \Bl(N_G(D)|D)$  be its Brauer correspondent block.
Then 
 $\k_{rv}(B)=\k_{rv}(B_0)$. 
\end{lemma}

\noindent
Proof. We use the same notation as in  the introduction to this section.
 By  \cite[Lemma 1.1]{Mur}  and Lemma \ref{cyclicdefect} it is enough 
to prove that if $(x^{2^{j}k},\theta^j)$  represents  a real column of the block $B$,
(recall that $\theta^j\in \IBr(b_j)$ and $b_j\in \Bl(C_G(D_j)|D)$),  
  and  $\tilde{b_j}\in \Bl(N_{C_G(D_j)}(D)|D)=\Bl(C_{N_G(D)}(D_j)|D)$ is the 
Brauer correspondent  of $b_j$    containing the single irreducible Brauer character
 $\tilde\theta^j$, then 
$(x^{{2^j}k},\tilde \theta^j)$  belongs to a real column of $B_0$ and  this correspondence defines a bijection
of real columns of $B_0$ and $B$.
 
Let $z\in G$ such that $(({x^{{2^j}k}})^z,{\theta^j}^z)=(({x^{{2^j}k}})^{-1},
\overline{\theta^j})$.
Then $\overline{\theta^j}\in \IBr(\overline{b_j})$.  This
block's  Brauer correspondent in $N_{C_G(D_j)}(D)$ is 
$\overline{\tilde{b_j}}$, that contains the unique
 irreducible Brauer character $\overline{\tilde{\theta^j}}$.
Since blocks of $C_{N_G(D)}(D_j)$ that induce $B_0$ are conjugate in $N_{N_G(D)}(D_j)=
N_G(D)$, there exists an element $z_1\in N_G(D)$ with  ${\tilde{b_j}}^{z_1}=\overline
{\tilde{b_j}}$. Then ${b_j}^{z_1}=\overline {b_j}$
 and ${\tilde {\theta^j}}^{z_1}=\overline{\tilde{\theta^j}}$.
 But then ${b_j}^{zz_1}={b_j}$. But the inertia group  of 
$b_j$ in $N_j$ is $C_j$, thus $zz_1\in C_j$ and so $({x^{{2^j}k}})^{zz_1}=x^{{2^j}k}$, and
 $({x^{{2^j}k}})^{-1}=({x^{{2^j}k}})^z=({x^{{2^j}k}})^{z_1}$, and hence
 $(({x^{{2^j}k}})^{z_1},{\tilde {\theta^j}}^{z_1})=(({x^{2^{j}k}})^{-1},\overline{\tilde{\theta^j}})$. Thus it represents a real column of $B_0$. 
The remaining column of $B$   containing  $(1,\phi)$   is real and the corresponding column
 containing 
$(1,\tilde\phi)$ in $B_0$ is also real. So we are done.

Now we have:

\begin{theorem}
Let $G$ be a finite group, let $B\in \Bl(G|D)$ be a $2$-block with cyclic defect group
$D$. Then Conjecture \ref{conj1} and hence Conjecture \ref{conj3} also holds for $B$.
\end{theorem}

\noindent
Proof.
  Using   Lemma \ref{samenumber} and  Corollary \ref{normal}, we have that
 $|\Irr_{rv}(B)|=|\Irr_{rv}(B_0)|$ is bounded  from above by the number of
the $N_G(D)$-real elements of $D$, thus  
 we are done.

\section{Groups with odd real core}

In \cite{HH} we defined the {\bf real core} $R(G)$ of $G$ as the subgroup
generated by the real elements of odd order.

\bigskip
Our main result is the following:

\begin{theorem}\label{theorem1}
If $|R(G)|$ is odd then Conjecture~\ref{conj3} holds (hence also Conjectures
\ref{conj2}, \ref{conj1}
and \ref{conj0})  for the principal $2$-block of
 $G$.  
In particular if any of the following cases occur Conjecture~\ref{conj3} holds for
the principal $2$-block of $G$. 
\begin{itemize}
\item [(a)] The commutator subgroup  $G'$ is $2$-nilpotent.
\item [(b)] $G=O_{2',2,2'}(G)$. (In fact this is equivalent to $|R(G)|$ being odd.)
\item [(c)] $G$ is solvable and its Sylow $2$-subgroup is abelian.
\end{itemize} 
Moreover, we may replace in Conjecture~\ref{conj3} $N_G(D)$ by $D$.
\end{theorem}

\noindent
Proof.  If $|R(G)|$ is odd then by \cite{HH} $G=O_{2',2,2'}(G)$, in particular
 $G$ is solvable.
Let $B_0$ be the principal $2$-block of $G$.
  Then $\Irr_{rv}(B_0)=\Irr_{rv}(G/O_{2'}(G))$. 
Let $\overline G=G/O_{2'}(G)$. Then $\overline S\in Syl_2(\overline G)$ is normal.

\noindent
Step 1 : For every  real element $x$ there exists a $2$-element $g$
 such that $x^g=x^{-1}$:

Let $g=g_2g_{2'}$ if $x^{g_2g_{2'}}=x^{-1}$, then an appropriate $2$-power
of $g$ is already a $2'$-element, which centralizes $x$. Thus $g_{2'}$ acts on $x$
trivially, and  we are done.

\noindent
Step 2: $R(\overline G)=1$:

If $x\in \overline G$ is a real element of odd order, then by Step 1 there is
a $2$-element g inverting $x$. Since $\overline S\triangleleft \overline G$,
  $[x,g]\in \overline S\cap \langle x \rangle =1$.
Thus $x^{-1}=x$,  and so $x=1$.

\noindent
Step 3: Every real element in $\overline G$ is a $2$-element, hence it lies in
$\overline S\in Syl_2(\overline G)$:

Let $x=x_2x_{2'}$ be a real element in $\overline G$. Then ${x_2}^{-1}{x_{2'}}^{-1}=
x^{-1}=x^g={x_2}^g{x_{2'}}^g$, thus $x_2$ and $x_{2'}$ are both real. By Step 2,
$x_{2'}=1$.

Thus $|\Irr_{rv}(B_0)|=|\Irr_{rv}(\overline G)|=|\Cl_{r}(\overline G)|\leq |\{
x\in\overline G | x$ ${\rm real} \}|$$=|\{ x \in \overline S\in Syl_2(\overline G)|
x $ ${\rm real}$ ${\rm in}$ $\overline S \}$.

We prove Conjecture~\ref{conj3} by induction on $r$. Let $r=0$. An irreducible character $\chi \in \Irr(B_0)$ is of height zero iff its degree is odd.
 We have that $\chi \in \Irr(\overline G)$, and $\chi_{\overline S}$ has only linear constituents, hence
${\overline{S}}'\leq \ker(\chi)$. Thus $\k_{0,rv}(B_0)\leq |\Irr_{rv}(\overline G/{\overline {S}}')|=|\Cl_r(\overline G/{\overline {S}}')|\leq |\{ x \in \overline G/{\overline {S}}'| x$ ${\rm real} $ ${\rm in}$ $\overline G/{\overline {S}}'$ $ \}|$.
If there would be a real $2'$-element in $G/{\overline {S}}'$ then by
 Proposition 5.3 in \cite{HH} there would be a real $2'$ element in $\overline G\backslash
{\overline {S}}'$, which is not the case. 
Thus there are also no real elements in $\overline G/{\overline {S}}'$ whose
 $2'$-part is not $1$. Thus each real element belongs to
 $\overline S/{\overline {S}}'$, and by Step 1 this element is 
 $\overline S /{\overline {S}}'$-real. 
 Thus we are done for $r=0$.

Let us suppose that Conjecture~\ref{conj3} is true for $r<n$.
If $\chi \in \Irr(B_0)$ is of height $n$, then its degree has $2$-part $2^n$.
Then all constituents of $\chi_{\overline S}$ have  degree $2^n$. By  
\cite[Thm.~5.12]{Is}, they contain in their kernels $\overline {S}^{(n+1)}$,
 thus $\chi$ also contains it in its kernel. Similarly all irreducible characters of
$\overline G$ of smaller height also contain it in their kernels.
Hence $\sum_{i=0}^n k_{i,rv}(B_0)\leq  |\Irr_{rv}(\overline G/{\overline S^{(n+1)})|=
|\Cl_r(\overline G/{\overline S^{(n+1)}}|\leq |\{x\in \overline G/{\overline S}^{(n+1)}}| x$  {\rm real}$\}|=
| \{x\in \overline S/{\overline S}^{(n+1)}| x $ {\rm real in\ }$ 
 {\overline S}/{\overline S}^{(n+1)}\}|$. Hence
Conjecture~\ref{conj3} holds.

\begin{itemize}
\item[(a)] Since $R(G)\leq G'$ by \cite{HH}, if $G'$ is $2$-nilpotent, then $|R(G)|$ is odd.
\item[(b)] This is equivalent to $|R(G)|$ odd by \cite{HH}.
\item[(c)] By the Hall-Higman lemma, $\overline S$ is normal, thus we have  case (b).
\end{itemize}

\noindent

\begin{corollary}
If $S\in Syl_2(G)$ is normal then Conjecture $\ref{conj3}$ holds for $G$, since
then each block is of maximal defect, and the only real $2$-block of maximal defect is the principal block, hence  we can apply Theorem \ref{theorem1} (b). 
\end{corollary}

\section{Computer results}

\bigskip

We have checked  Conjecture~\ref{conj3} for the principal $2$-block 
 with \textsf{GAP} \cite{GAP} for the small groups library. We   also checked
Conjecture~\ref{conj1} for the principal $2$-block 
for the $26$ sporadic simple groups.
For these blocks the respective conjectures were true.


We also checked Conjecture \ref{conj1} 
 and  Conjecture \ref{conj2} with the help of $GAP$ for all
$2$-blocks of groups up to order $1536$ except for  groups of orders
 $856,1048,1112,1192,1304,1352,1384,1432,1448$, where our computational methods did not work (Conway polynomials are not yet known).
We did not find any conterexamples for these conjectures among the investigated groups.

\section{$B$-classes  and $G$-classes of $D$}

Let now $p$ be an arbitrary prime number.

First we prove the following

\begin{lemma}\label{Hiss}
Let $B$ be a $p$-block of $G$ with defect group $D$.
Then  for every $x\in D$  there exists a 
 $\chi \in \Irr(B)$ with $\chi(x)\ne 0$.
\end{lemma}

\noindent
Proof. Let us suppose by contradiction that there exists
an element $x\in D$ with $\chi(x)=0$, for every $\chi\in \Irr(B)$. 
If we can prove that                                                                   there exists a trivial source $FG$-module 
$M$ in this block with vertex $D$,
 then by \cite[p.~175, Lemma~2.16]{Lan} this is 
liftable to a trivial source $RG$-module $\tilde M$ and  its                            character is nonzero on the elements of  the vertex of $M$,
contradicting  our assumption.

If $D$ is normal in G then by \cite[p.~247, Lemma~10.3]{Lan} 
all simple modules in
$B$ are trivial source modules in $B$ with vertex $D$.
If $D$ is not normal then  the Brauer correspondent
$b\in \Bl(N_G(D)|D)$ of $B$ has the property that every simple module
 $S$ in it is a trivial source module with  vertex $D$.
Let us lift a simple $FN_G(D)$-module $S$ in $b$ to an $RN_G(D)$-module
$\tilde S$. Then its Green correspondent, $f(\tilde S)$ is a
 trivial source module
with
vertex $D$ and by \cite[p.~466, Thm.~59.9]{C-R}, $f(\tilde S)$ belongs 
to the block $B$.
So we are done.
 
\noindent
\begin{definition}  Let $B$ be a $p$-block of the finite group $G$ with defect group 
$D$. We say that two elements $x,y\in D$ are in the same $B$-class, iff for every
irreducible character $\chi\in \Irr(B)$, $\chi (x)=\chi (y)$.
\end{definition}

 We have the following result:

\begin{theorem}\label{bclass} Let $G$ be a finite group with $p$-block $B\in \Bl(G|D)$.
Then the $B$-classes of the defect group $D$ are exactly the $G$-classes $\Cl_G(D)$
of $D$ under conjugation.
\end{theorem}

\noindent
Proof. If two elements $x,y\in D$ are conjugate in $G$, then of course they are also in the same $B$-class.
Let us suppose now that $x,y\in D$ are in the same $B$-class, but they
 are not conjugate in $G$. Then by the strong
block orthogonality relation, see \cite[p.~106, Cor.~5.11]{Nav} $\sum_{\chi\in \Irr(B)}\chi(x)\overline{\chi(y)}=0$.
Using that $x,y$ are in the same $B$-class this gives us
$\sum_{\chi(x)\in \Irr(B)}|\chi(x)|^2=0$.
Hence $\chi(x)=0$ for every $\chi\in \Irr(B)$. 
By Lemma \ref{Hiss}, this is not possible.

\begin{definition} Let $B$ be a $p$-block of a finite group
 $G$ with defect group $D$. We say that the element $x\in D$ is
 $B$-real if $\chi(x)$ is real for every $\chi\in \Irr(B)$.
 An element $x\in D$ is $B$-rational, if $\chi(x)$ is rational 
for every
$\chi \in \Irr(B)$. 
\end{definition}

\begin{corollary} Let $B\in \Bl(G|D)$ and let $x\in D$. Then $x$ is
 $B$-real iff
it is real in $G$.
\end{corollary}

\begin{corollary} Let $B\in \Bl(G|D)$ and let $x\in D$. Then $x$ is $B$-rational
iff it is rational in $G$.
\end{corollary}

\begin{corollary} Let 
Let $F$ be  a field containing $Q$.
Let $B\in Bl(G|D)$ and let $x\in D$. Then $\chi (x)\in F$ for every $\chi \in \Irr(G)$
if and only if  $\chi (x)\in F$ for every $\chi \in \Irr(B)$.
\end{corollary}

We have also the following

\begin{theorem}
 Let $B\in \Bl(G|D)$. The restrictions $\chi_D$ of $\chi\in \Irr(B)$
 to the defect group $D$ generate the vector space of the restrictions
of all complex $G$-class 
functions to $D$.
\end{theorem}

\noindent
Proof.
Let us choose representatives $x_i\in C_i\cap D, i = 1, \ldots, t$
 of $G$-conjugacy classes $C_i$ intersecting $D$.
We want to prove that if we restrict the character table of $G$ to
 these columns
and to those rows which belong to the block $B$, then these columns
are independent.  It implies that this
 submatrix has  rank $t$, hence, it has also
$t$  independent rows. But then any complex vector 
of length $t$ can be expressed by these rows and we are done.
Let us suppose that the above mentioned $t$ columns are dependent.
 Then
there are coefficients $\alpha_1,...,\alpha _t$, not all zero 
with the property
that
$\sum_{i=1}^t\alpha_i\chi(x_i)=0$, for all $\chi \in \Irr(B)$.
By \cite[Lemma~4.6, Ch.~5]{N-T}   the subsum, where $x_i$-s belong to
any
$p$-section  is also zero. But the $x_i$-s all belong to different 
$p$-sections,
thus $\alpha_i\chi(x_i)=0$ for every $i = 1, \ldots, t$ and 
every $\chi \in \Irr(B)$.
By Lemma \ref{Hiss}  we see that there is no $x_i$ where
every
$\chi \in \Irr(B)$ vanishes. Hence $\alpha _i=0$  for all
 $i = 1, \ldots, t$. Thus
the columns of the above restricted matrix are independent and
 we are done.

In this way we get another proof of the following:

\begin{corollary}\label{lowerbound} For a block $B\in \Bl(G|D)$, the number of $G$-conjugacy classes 
$|\Cl_G(D)|$ of its defect group, is a lower bound for the number $\k(B)$.
\end{corollary}



\begin{example}
It is not true, however that $B$-classes (hence $G$-classes) of the defect group $D$
 are the same as $b$-classes (hence $N_G(D)$-classes) 
for the Brauer correspondent block $b\in \Bl(N_G(D)|D)$ even for
$2$-blocks $B\in \Bl(G|D)$ with cyclic defect group.
Let $G=\texttt{SmallGroup}(288,375)$. Then the third $2$-block has cyclic defect group of order $8$,
it contains four $G$-real (hence $B$-real)  elements  and only  two $N_G(D)$-real
(hence $b$-real) elements.

\end{example}

{\rm ACKNOWLEDGEMENTS}
The authors are very grateful to Professor Gerhard Hiss for  the
discussions about the topic of this paper and also
 for his suggestions and comments, especially for his
 contributions to the proof of Lemma \ref{Hiss}.  
We also thank Thomas Breuer for his advices concering the GAP system. 
Research supported by National Science Foundation Research Grant No. T049841
and K 77476.

\bibliographystyle{amsalpha}

\begin{thebibliography}{ABCD}


\bibitem{Brau} R.~Brauer:
  \newblock {\em Some applications of the theory of blocks and characters
            of finite groups IV.},
  \newblock {J. Algebra} {\bf 17} {\bf 4} (1971), 20--52.

\bibitem{Bra} R.~Brauer:
  \newblock  Representations of finite groups, Lectures on Modern Mathematics
  (ed. by T. Saaty) Vol 1. pp.~133--175. Wiley, New York, 1963.

\bibitem{Bra56} R.~Brauer:
\newblock {\em Number theoretical investigations on groups of finite order},
\newblock {Proc. Int. Symp. Alg. Number Theory}, Tokyo Nikko, (1956), 55--62.

\bibitem{C-R} C.W. Curtis, I. Reiner:
\newblock Methods of representation theory, Vol II., Wiley and Sons, New York, 1987. 

\bibitem{DolNavTie} S.~Dolfi, G.~Navarro and P.~H.~Tiep:
 \newblock {\em Primes dividing the degrees of the real characters},
  \newblock{Math.Z.} {\bf 259} (2008), 755--744.

\bibitem{Dor} L.~Dornhoff, Group representation theory Part B,
\newblock Marcel Dekker Inc.,  New York 1972.

\bibitem{GAP}  The \textsf{GAP} Group,
   \newblock \textsf{GAP}--Groups, Algorithms and Programming,
   
\newblock Version 4.4.12, \texttt{http://www.gap-system.org}, 2008.

\bibitem{Gow} R.~Gow:
  \newblock {\em Real $2$-blocks of characters of finite groups},
  \newblock {Osaka J. Math.} {\bf 25} (1988), 135--147.

\bibitem{G-M} R.~Gow and J.~Murray:
  \newblock {\em Real $2$-regular classes and $2$-blocks},
  \newblock {J. Algebra} {\bf 230} (2000), 455--473.

\bibitem{Eat03} C.~W.~Eaton:
\newblock {\em Generalisations of conjectures of Brauer and Olsson},
\newblock
{Arch. Math.} {\bf 81} (2003), 621--626.


\bibitem{HH} L.~H\'ethelyi and E.~Horv\'ath:
\newblock {\em Galois actions on blocks and classes of finite groups},
\newblock{J. Algebra} {\bf 320} (2008), 660--679.

\bibitem{Is} I.~M.~Isaacs:
 \newblock 
Character theory of finite groups,
\newblock Academic Press, New York, 1976.

\bibitem{IMN} I.~M.~Isaacs, G.~Malle and G.~Navarro:
\newblock {\em Real characters of $p'$ degree},
\newblock{J. Algebra} {\bf 278} (2004), 611--620.

\bibitem{IN} I.~M.~Isaacs, G.~Navarro:
\newblock {\em Group elements and fields of character values},
\newblock{J. Group Theory} {\bf 12(5)} (2009), 625--650.

\bibitem{Iwasaki} S.~Iwasaki:
\newblock {\em On finite groups with exactly two real conjugacy classes},
\newblock { Arch. Math.} {\bf 33} (1979), 512--517.

\bibitem{J-K} G.~D.~James, A.~Kerber, The Representation Theory of
 the Symmetric Group,
\newblock Addison-Wesley Publ. Co., London, 1981.


\bibitem{K-N} S.~G.~Kolesnikov and J.~A.~N.~Nuzhin:
  \newblock {\em On strong reality of finite simple groups},
  \newblock {Acta Applicandae Mathematicae} {\bf 85} (2005), 195--203.

\bibitem{Lan} P.~Landrock:
\newblock Finite group algebras and their modules,
\newblock  LMS Lecture Note Series 84, CUP, London, 1983.

\bibitem{MN} A.~Moreto and G.~Navarro:
\newblock {\em Groups  with three  real valued irreducible characters}, 
\newblock {Israel J.  Math.} {\bf 163} (2008), 85--92.

\bibitem{M} J.~Murray:
  \newblock {\em Involutions, extended defect groups of real $2$-blocks
  and vertex theory},  Preprint.


\bibitem{Mur} J.~Murray:
\newblock  {\em Real subpairs and Frobenius-Schur indicators of characters
in $2$-blocks}, 
\newblock {J. Algebra} {\bf 322(2)} (2009), 489--513.

\bibitem{N-T} H.~Nagao and Y.~Tsushima:
 \newblock Representations of finite groups,
 \newblock Academic Press, New York, 1988.

\bibitem{Nav} G.~Navarro:
  \newblock Characters and Blocks of Finite Groups,
  \newblock LMS Lecture Note Series, Vol 250, Cambridge University Press,
   Cambridge, 1998.

\bibitem{NavSanTie} G.~Navarro, L.~Sanus, P.~H.~Tiep:
\newblock {\em Groups with two real Brauer characters},
\newblock {J. Algebra} {\bf 307} (2007), 891--898.


\bibitem{NavSanTie2} G.~Navarro, L.~Sanus, P.~H.~Tiep:
\newblock {\em Real characters and degrees},
\newblock {Israel J. Math.} to appear.


\bibitem{Ol76} J.~Olsson:
\newblock {\em McKay numbers and heights of characters},
\newblock {Math. Scand.} {\bf 38} (1976), 25--42.

\bibitem{Ol80} J.~Olsson:
\newblock {\em Inequalities for block theoretic invariants},
\newblock Representations of Algebras Proceedings, Puebla, Mexico, 1980,
\newblock Lecture Notes in Mathematics 903, Ed. A.~Dold and  B.~Eckmann,
\newblock Springer-Verlag, Berlin, Heidelberg, New York, 1981, 270--284.  

\bibitem{Ol84} J.~Olsson:
\newblock {\em On the number of characters in blocks of finite general linear, unitary and symmetric groups},
\newblock {Math. Z.} {\bf 186} (1984), 41--47.

\end{thebibliography}

\end{document}